\documentclass[11pt]{article}

\usepackage{amsmath,amssymb,amsthm,esint,mathtools}
\usepackage{amssymb}

\usepackage{color,enumitem,graphicx,cite,mathrsfs,tikz,colortbl}

\usepackage[colorlinks=true,urlcolor=blue,citecolor=red,linkcolor=blue,linktocpage,pdfpagelabels,bookmarksnumbered,bookmarksopen]{hyperref}
\usepackage[english]{babel}

\usepackage[pdftex,a4paper,left=2.9cm,right=2.9cm,top=3.2cm,bottom=3.2cm]{geometry}

\usepackage[hyperpageref]{backref}
\usepackage{cite}
\usepackage{setspace}
\newtheorem{theorem1}{Theorem}[section]
\newtheorem{lemma1}[theorem1]{Lemma}

\theoremstyle{remark}

\numberwithin{equation}{section}

\title{\bf Existence result for fractional problems with logarithmic and critical exponential nonlinearities \thanks{{\em MSC2020:} Primary 35J35 Secondary 35J60, 35D30
\newline \indent\; 
{\em Key Words and Phrases:} Fractional Trudinger-Moser embedding, Logarithmic nonlinearity, Exponential nonlinearity, Existence of solutions}}

\author{\bf Yuanyuan Zhang\\
School of  Science\\
Jiangnan university\\
Wuxi, Jiangsu  214122, China\\
E-mail: 6191204010@stu.jiangnan.edu.cn\\
[\bigskipamount]
\bf Yang Yang\thanks{Corresponding author. Project supported by NSFC(No. 11501252, No. 11571176).}\\
School of Science\\
Jiangnan University\\
Wuxi, Jiangsu  214122, China\\
E-mail: yynjnu@126.com}
\date{}

\begin{document}
\begin{spacing}{1.25}
\maketitle
\begin{abstract}
    We study the existence of nontrivial solutions for a nonlinear fractional elliptic equation in presence of 
    logarithmic and critical exponential nonlinearities. This problem extends \cite{Liang2021} to fractional
    $p$-Laplacian equations with logarithmic nonlinearity.
    We overcome the lack of compactness due to the critical exponential nonlinearity by using
    the fractional Trudinger-Moser inequality. The existence result is established via critical point theory.
\end{abstract}


\section{Introduction and main results}

Let $N\geq 2$, $s\in (0,1)$ and $\Omega \subset \mathbb{R}^N$ be a bounded domain with Lipschitz boundary.
For $\mu>0$ and $\frac{N}{s}< q < \frac{N\sigma}{s}$, we consider the existence of solutions 
for the following fractional problem
\begin{equation}
    \label{main equations}
    \left\{
        \begin{array}{llll}
            (-\Delta)_{N/s}^s u = |u|^{q-2}u\ln |u|^2 + \mu |u|^{(N\sigma - 2s)/s}u e^{|u|^{N/(N-s)}} , & \text{in}\,\Omega,\\
            u=0, & \text{in}\,\mathbb{R}^N\backslash \Omega.
        \end{array}
    \right.
\end{equation}
Here, $(-\Delta)_{N/s}^s$ is the nonlinear nonlocal operator defined on smooth functions by
\begin{equation*}
    (-\Delta)_{N/s}^s \phi (x) \coloneqq 2 \lim_{\epsilon\searrow 0} \int_{\mathbb{R}^N \backslash B_{\epsilon}(x)}
    \frac{|\phi (x)-\phi (y)|^{\frac{N}{s}-2}(\phi (x)- \phi (y))}{|x-y|^{2N}} \,dy, \,\, x\in \mathbb{R}^N,
\end{equation*}
for all $\phi \in \mathcal{C}_0^\infty(\Omega)$, where $B_{\epsilon}(x) \subset \mathbb{R}^N$ denotes an open ball centered at $x\in \mathbb{R}^N$ with
radius $\epsilon > 0$.

Since the fractional Trudinger-Moser inequality (see Lemma \ref{fractional Trudinger-Moser inequality}) 
was established, the existence and multiplicity of solutions for various nonlinear nonlocal problems with exponential 
nonlinearity were investigated.
We refer to \cite{Xiang-Zhang2019,Chen-Yu2020,Miyagaki-Pucci2019,Souzaa-Araujo2017,Souza-Severo2019,Souza-Severo2021} 
and references therein.
In particular, Perera and Squassina \cite{Perera-Marco2018} obtained the bifurcation results of the following 
problems 
\begin{equation*}
    \left\{
        \begin{array}{llll}
            (-\Delta)_{N/s}^s u = \lambda |u|^{(N-2s)/s} u \exp (|u|^{N/(N-s)}) , & \text{in}\,\Omega,\\
            u=0, & \text{in}\,\mathbb{R}^N\backslash \Omega,
        \end{array}
    \right.
\end{equation*}
for appropriate $\lambda>0$.
On the other hand, there have been several papers that deal with fractional problems involving logarithmic nonlinearities. 
Xiang, Hu and Yang \cite{Xiang-Hu-Yang2020} considered the following Kirchhoff problems:
\begin{equation*}
    \left\{
        \begin{array}{llll}
            M([u]_{s,p}^p)(-\Delta)_p^s u = h(x)|u|^{\theta p -2} u \ln |u| + |u|^{q-2}u , & \text{in}\,\Omega,\\
            u=0, & \text{in}\,\mathbb{R}^N\backslash \Omega,
        \end{array}
    \right.
\end{equation*}
for $1 < p < N/s$, $\theta \in (1, p^*_s/p)$, $q\in(1,\theta p)$ and appropriate $h \in C(\overline{\Omega})$,
where $p^*_s = \frac{Np}{N-sp}$ is the fractional critical Sobolev exponent.
Liang et al \cite{Liang2021}, subsequently, studied the case $q=p^*_s$.
However, to our best knowledge, the existence of solutions for fractional $N/s$-Laplacian equations with logarithmic and 
critical exponential nonlinearities was not tackled in the literature, which is precisely the main goal of this manuscript.

Here is our main existence result.
\begin{theorem1}
    \label{main theorem}
    Assume that $N\geq 2$, $ 0<s<1 $ and $\frac{N}{s}< q < \frac{N \sigma}{s}$. 
    Then there exists $\mu^{*}>0$ such that for any $\mu \geq \mu^*$ problem (\ref{main equations}) admits a nontrivial solution.
\end{theorem1}

The outline of the paper is as follows:
Section 2 contains the variational framework and some auxiliary results. 
And we also show that the corresponding energy functional satisfies the Palais-Smale condition below appropriate critical level.
In Section 3, we prove the existence of at least one solution for problem (\ref{main equations}) 
using some asymptotic estimates and the Mountain Pass Theorem.

\section{Auxiliary results} 
We start this section with the fractional Sobolev space and some information about the weak formulation of problem (\ref{main equations}).
Throughout this paper, we work in the fractional Sobolev space 
$W_0^{s,N/s}(\Omega)$, defined as the completion of $\mathcal{C}_0^\infty(\Omega)$ 
with respect to the Gagliardo seminorm
\begin{equation*}
    \|u\| = [u]_{s,N/s} = (\int_{\mathbb{R}^{2N}} \frac{|u(x)-u(y)|^{N/s}}{|x-y|^{2N}} \,dxdy)^{s/N}.
\end{equation*}

The well-known Sobolev embedding theorem states that for $\nu \in[1,\infty)$, 
$W_0^{s,N/s}(\Omega) \hookrightarrow L^{\nu}(\Omega)$ is continuous and compact. One can refer to 
\cite{Hitchhiker2012} for more detailed account of the properties of $W_0^{s,N/s}(\Omega)$.
To study problems involving exponential growth in the fractional Sobolev space, the main tool is the fractional
Trudinger-Moser inequality due to \cite{Parini-Ruf2019} (see also \cite{Perera-Marco2018,ZhangCaifeng2019}).

\begin{lemma1}{\cite[Theorem 1.1]{Parini-Ruf2019}}
    \label{fractional Trudinger-Moser inequality}
    Let $\Omega$ be a bounded, open domain of $\mathbb{R}^N$($N \geq 2$) with Lipschitz boundary, and let $s\in(0,1)$.
    Then there exists $\alpha_{s,N}$ such that
    \begin{equation*}
        \sup \{ \int_{\Omega}\exp (\alpha |u|^{\frac{N}{N-s}}) \,dx : u \in W_0^{s,N/s}(\Omega), [u]_{s,N/s} \leq 1
         \} < \infty,
    \end{equation*}
    for $\alpha \in [0,\alpha_{s,N})$. Moreover,
    \begin{equation*}
        \sup \{ \int_{\Omega}\exp (\alpha |u|^{\frac{N}{N-s}}) \,dx : u \in W_0^{s,N/s}(\Omega), [u]_{s,N/s} \leq 1
         \} = \infty,
    \end{equation*}
    for $\alpha \in (\alpha^*_{s,N}, \infty)$, 
    where
    \begin{equation*}
        \alpha^*_{s,N} \coloneqq N (\frac{2(N\omega_N)^2 \Gamma (N/s+1)}{N!} 
        \sum_{k = 0}^{\infty} \frac{(N+k-1)!}{k!} \frac{1}{(N+2k)^{N/s}})^{\frac{s}{N-s}}.
    \end{equation*}
    and $\omega_N$ is the volume of the $N$-dimential unit ball.
\end{lemma1}

\begin{lemma1}{\cite[Proposition 3.2]{Parini-Ruf2019}}
    \label{Corollary of fractional Trudinger-Moser inequality}
    If $u \in W_0^{s,N/s}(\Omega)$, then for any $\alpha \geq 0$, it holds
    \begin{equation*}
        \int_{\Omega} \exp (\alpha |u|^{\frac{N}{N-s}}) \,dx < \infty.
    \end{equation*}
\end{lemma1}

Let $f(t) = |t|^{(N\sigma - 2s)/s} t e^{|t|^{N/(N-s)}}$ for $t\in \mathbb{R}$ and $F(t)= \int_{0}^{t} f(\tau ) \,d\tau$.
For each $\mu >0$, a function $u \in W_0^{s,N/s}(\Omega)$ is a weak solution of problem (\ref{main equations}) if 
\begin{equation*}
    \begin{split}
        \int_{\mathbb{R}^{2N}} & \frac{|u(x)-u(y)|^{(N-2s)/s}(u(x)-u(y))(v(x)-v(y))}{|x-y|^{2N}}\,dxdy\\
        &= \int_{\Omega}|u|^{q-2}uv \ln |u|^2 \,dx + \mu \int_{\Omega}f(u)v \,dx,\,\, \forall v \in W_0^{s,N/s}(\Omega).
    \end{split}
\end{equation*}
The last integral is well-defined in view of Lemma \ref{Corollary of fractional Trudinger-Moser inequality} and the
H$\ddot{\text{o}}$lder inequality. The weak solutions of problem (\ref{main equations}) coincide with the critical points of 
the $C^1$ functional 
\begin{equation*}
        I_{\mu}(u)= \frac{s}{N}\|u\|^{N/s} + \frac{2}{q^2} \int_{\Omega} |u|^q \,dx 
        - \frac{1}{q} \int_{\Omega} |u|^q \ln |u|^2 \,dx - \mu \int_{\Omega} F(u) \,dx, \,\, u \in W_0^{s,N/s}(\Omega).
\end{equation*}

It follows from $q > N/s$ that $\lim_{t \to 0} \frac{|t|^{q-1}\ln |t|^2}{|t|^{\theta_1 -1 }} = 0$ and 
$\lim_{t \to \infty} \frac{|t|^{q-1}\ln |t|^2}{|t|^{\theta_2 -1 }} = 0$ for all $\theta_1 \in (N/s, q)$
and $\theta_2 \in (q,\infty)$. Then for any $\epsilon >0 $, there exists $C_{\epsilon}>0$ such that
\begin{equation}
    \label{1,estimate of the logarithmic nonlinearity}
    |t|^{q-1} \ln |t|^2 \leq \epsilon |t|^{\theta_1 -1 } + C_{\epsilon} |t|^{\theta_2 -1}.
\end{equation}
Besides, we observes that 
\begin{equation}
    \label{2,estimate of the logarithmic nonlinearity}
    2t^q -qt^q \ln |t|^2 \leq 2, \,\, \text{for all}\,\, t\in (0,\infty),
\end{equation}
by simple calculation.

Next, we will collect some elementary estimates of the exponential nonlinearity.
\begin{lemma1}
    \label{estimate of the exponential nonlinearity}
    For all $t\in \mathbb{R}$,
    \begin{itemize}
        \item [(i)] $F(t) \leq \frac{1}{q} f(t)t$;
        \item [(ii)] $F(t) \leq \frac{s}{N \sigma} |t| ^{\frac{N \sigma}{s}} + |t|^{\theta} e^{|t|^{N/(N-s)}}$;
        \item [(iii)] $F(t) \geq \frac{s}{N \sigma} |t| ^{\frac{N \sigma}{s}} + \frac{1}{\theta}|t|^{\theta}$,
        where $\theta = \frac{N \sigma}{s} + \frac{N}{N-s}$.
    \end{itemize}
\end{lemma1}

\begin{proof}
    Since $f$ is odd, $F(t) = \int_0^t f(\tau)\,d\tau = \int_0^{|t|}  \tau ^{(N \sigma -s)/s} e^{\tau^{N/(N-s)}}\,d\tau$.

    (i) Integrating by parts,
    \begin{equation*}
        \begin{split}
            F(t)
            & = \frac{s}{N \sigma} tf(t) - \frac{s}{\sigma (N-s)} \int_0^{|t|}  \tau ^{\theta-1} e^{\tau^{N/(N-s)}}\,d\tau 
            \leq \frac{1}{q} tf(t).
        \end{split}
    \end{equation*}

    (ii) Since $e^{\tau} \leq 1 + \tau e^{\tau}$ for all $\tau \geq 0$,
    \begin{equation*}
        \begin{split}
            F(t)
            & \leq \int_0^{|t|}  \tau ^{(N \sigma -s)/s} \,d\tau + 
            \int_0^{|t|} |t|^{\theta -1} e^{|t|^{N/(N-s)}}  \,d\tau
            \leq \frac{s}{N \sigma} |t| ^{\frac{N \sigma}{s}} + |t|^{\theta} e^{|t|^{N/(N-s)}}.
        \end{split}
    \end{equation*}

    (iii) Since $e^{\tau} \geq 1 + \tau$ for all $\tau \geq 0$,
    \begin{equation*}
        \begin{split}
            F(t)
            & \geq \int_0^{|t|}  \tau ^{(N \sigma -s)/s} (1+\tau^{N/(N-s)})\,d\tau  
            = \frac{s}{N \sigma} |t| ^{\frac{N \sigma}{s}} + \frac{1}{\theta}|t|^{\theta}.
        \end{split}
    \end{equation*}
\end{proof}

The next result is concerned in the $\text{(PS)}_c$ condition.

\begin{lemma1}
    \label{PS_c condition}
    $I_{\mu}$ satisfies the $\text{(PS)}_c$ condition for all $0< c < (\frac{s}{N} - \frac{1}{q}) \alpha_{s,N}^{(N-s)/s}$.
\end{lemma1}

\begin{proof}
    Let $(u_j)\subset W_0^{s,N/s}(\Omega)$ satisfy $I_{\mu} (u_j) \to c$ and $I'_{\mu}(u_j)\to 0$ for 
    $0< c < (\frac{s}{N} - \frac{1}{q}) \alpha_{s,N}^{(N-s)/s}$. If $\inf\limits_j \|u_j\| = 0$, then there exists a subsequence of $(u_j)$, 
    still denoted by $(u_j)$, such that $u_j \to 0$ in $W_0^{s,N/s}(\Omega)$ as $j \to \infty$, which is impossible for $c>0$.
    Thus we assume that $\inf\limits_j \|u_j\| > 0$.

    It follows from Lemma \ref{estimate of the exponential nonlinearity} (i) that
    \begin{equation}
        \label{PS sequence is bdd}
        \begin{split}
            c+o(1)\|u_j\| 
            & = I_{\mu}(u_j) - \frac{1}{q}\langle I'_{\mu}(u_j), u_j \rangle \\
            & = (\frac{s}{N}-\frac{1}{q})\|u_j\|^{N/s} + \frac{2}{q^2} \int_{\Omega} |u_j|^q \,dx 
            + \mu \int_{\Omega} [\frac{1}{q}f(u_j)u_j-F(u_j)] \,dx \\
            & \geq (\frac{s}{N}-\frac{1}{q})\|u_j\|^{N/s}.
        \end{split}
    \end{equation}
    By $\frac{N}{s}<q$, we deduce that the sequence $(u_j)$ is bounded in $W_0^{s,N/s}(\Omega)$.
    Therefore, there exists some 
    $u\in W_0^{s,N/s}(\Omega)$ such that $(u_j)$ converges to $u$ weakly in $W_0^{s,N/s}(\Omega)$, stronly in $L^{\nu}(\Omega)$ for $\nu \in [1,\infty)$
    and a.e. in $\Omega$.
    Then we deduce from (\ref{1,estimate of the logarithmic nonlinearity}), the Sobolev inequality and $ \langle I'_{\mu}(u_j), u_j \rangle \to 0$
    that $\int_{\Omega} |u_j|^q \ln |u_j|^2 \,dx$ and $\int_{\Omega} f(u_j)u_j \,dx$ are bounded.

    Since $I'_{\mu}(u_j)\to 0$,
    \begin{equation}
        \label{step 1 to I'_mu(u)=0}
        \begin{split}
            \int_{\mathbb{R}^{2N}} & \frac{|u_j(x)-u_j(y)|^{(N-2s)/s}(u_j(x)-u_j(y))(v(x)-v(y))}{|x-y|^{2N}}\,dxdy\\
            &-\int_{\Omega}|u_j|^{q-2}u_j v \ln |u_j|^2 \,dx - \mu \int_{\Omega}f(u_j)v \,dx \to 0,
        \end{split}
    \end{equation}
    for all $v\in W_0^{s,N/s}(\Omega)$. For $v \in \mathcal{C}_0^{\infty}(\Omega)$ and $M>0$, 
    \begin{equation*}
        \begin{split}
            \int_{\Omega} f(u_j)v \,dx
            & = \int_{\{|u_j|<M\}} f(u_j)v \,dx + \int_{\{|u_j|\geq M\}} f(u_j)v \,dx\\
            & \leq \int_{\{|u_j|<M\}} f(u_j)v \,dx + \frac{\sup |v|}{M}\int_{\Omega} f(u_j)u_j \,dx\\
            & = \int_{\{|u_j|<M\}} f(u_j)v \,dx + \mathcal{O}(\frac{1}{M}),
        \end{split}
    \end{equation*}
    which shows that 
    \begin{equation}
        \label{step 2 to I'_mu(u)=0}
        \int_{\Omega}f(u_j)v \,dx \to \int_{\Omega}f(u)v \,dx,
    \end{equation}
    by letting $j\to \infty$ first and then $M \to \infty$.
    Similarly, we also have that 
    \begin{equation}
        \label{step 3 to I'_mu(u)=0}
        \int_{\Omega}|u_j|^{q-2}u_j v \ln |u_j|^2 \,dx \to \int_{\Omega}|u|^{q-2}uv \ln |u|^2 \,dx,
    \end{equation}
    as $j \to \infty$.
    It follows from (\ref{step 1 to I'_mu(u)=0})-(\ref{step 3 to I'_mu(u)=0}) that
    \begin{equation*}
        \begin{split}
            \int_{\mathbb{R}^{2N}} & \frac{|u(x)-u(y)|^{(N-2s)/s}(u(x)-u(y))(v(x)-v(y))}{|x-y|^{2N}}\,dxdy\\
            &= \int_{\Omega}|u|^{q-2}u v \ln |u|^2 \,dx + \mu \int_{\Omega}f(u)v \,dx.
        \end{split}
    \end{equation*} 
    Then the last equality holds for all $v\in W_0^{s,N/s}(\Omega)$ by density and 
    \begin{equation}
        \label{u is critical point}
        \|u\|^{N/s} = \int_{\Omega}|u|^q \ln |u|^2 \,dx + \mu \int_{\Omega}f(u)u \,dx
    \end{equation}
    holds by taking $v=u$.

    Next, we claim that 
    \begin{equation}
        \label{f(u_j)u_j to f(u)u}
        \int_{\Omega} f(u_j)u_j \,dx \to \int_{\Omega}f(u)u \,dx.
    \end{equation}
    On one hand, by (\ref{PS sequence is bdd}), Lemma \ref{fractional Trudinger-Moser inequality}
    and choosing $t>1$ such that $t\|u_j\|^{N/(N-s)} < \alpha_{s,N}$, we can conclude that $\int_{\Omega} e^{t|u_j|^{N/(N-s)}} \,dx \leq M$.
    Since $ e^{t|u_j|^{N/(N-s)}} \to e^{t|u|^{N/(N-s)}}$ a.e. in $\Omega$, we use \cite[Lemma 4.8]{Kavian1991} and conclude that
    \begin{equation}
        \label{step 1 of proving f(u_j)u_j to f(u)u}
        e^{|u_j|^{N/(N-s)}} \rightharpoonup e^{|u|^{N/(N-s)}} \,\, \text{in} \,\, L^t (\Omega).
    \end{equation}
    On the other hand, note that
    \begin{equation}
        \label{step 2 of proving f(u_j)u_j to f(u)u}
        |u_j|^{N \sigma /s} \to |u|^{N \sigma /s} \,\, \text{in}\,\, L^{t'}(\Omega),
    \end{equation}
    for $t'>1$ with $\frac{1}{t}+ \frac{1}{t'}=1$.

    Then it follows from the H$\ddot{\text{o}}$lder inequality, (\ref{step 1 of proving f(u_j)u_j to f(u)u}),
    (\ref{step 2 of proving f(u_j)u_j to f(u)u}) and Lemma \ref{fractional Trudinger-Moser inequality} that 
    \begin{equation*}
        \begin{split}
            \int_{\Omega} [f(u_j)u_j-f(u)u] \,dx
            & = \int_{\Omega} [|u_j|^{N \sigma /s} - |u|^{N \sigma /s}]e^{|u_j|^{N/(N-s)}} \,dx\\
            & + \int_{\Omega} |u|^{N \sigma /s} [e^{|u_j|^{N/(N-s)}} - e^{|u|^{N/(N-s)}}] \,dx\\
            & \leq (\int_{\Omega} [|u_j|^{N \sigma /s} - |u|^{N \sigma /s}]^{t'} \,dx)^{1/t'} 
            (\int_{\Omega} e^{t|u_j|^{N/(N-s)}} \,dx)^{1/t} \\
            & + \int_{\Omega} |u|^{N \sigma /s} [e^{|u_j|^{N/(N-s)}} - e^{|u|^{N/(N-s)}}] \,dx\\
            & \leq M (\int_{\Omega} [|u_j|^{N \sigma /s} - |u|^{N \sigma /s}]^{t'} \,dx)^{1/t'} \\
            & + \int_{\Omega} |u|^{N \sigma /s} [e^{|u_j|^{N/(N-s)}} - e^{|u|^{N/(N-s)}}] \,dx\\
            & \to 0 \,\, \text{as}\,\, j\to \infty,
        \end{split}
    \end{equation*}
    which complete the proof of (\ref{f(u_j)u_j to f(u)u}).

    Now, we are ready to show that the sequence $(u_j)$ converges strongly to $u$ in $W_0^{s,N/s}(\Omega)$.
    Indeed, by the weak lower semicontinuity of the norm, the Lebesgue dominated convergence theorem, (\ref{1,estimate of the logarithmic nonlinearity}),
    (\ref{u is critical point}) and (\ref{f(u_j)u_j to f(u)u}),
    we obtain that
    \begin{equation*}
        \begin{split}
            0
            & = \lim_{j\to\infty} \langle I'_{\mu}(u_j), u_j\rangle \\
            & = \lim_{j\to\infty} [ \|u_j\|^{N/s} - \int_{\Omega} |u_j|^q \ln |u_j|^2 \,dx - \mu \int_{\Omega}f(u_j)u_j \,dx]\\
            & \geq \|u\|^{N/s} - \int_{\Omega} |u|^q \ln |u|^2 \,dx - \mu \int_{\Omega}f(u)u \,dx\\
            & = 0,
        \end{split}
    \end{equation*}
    then $(u_j)$ converges strongly to $u$ in $W_0^{s,N/s}(\Omega)$. This completes the proof.
\end{proof}

\section{Proof of Theorem \ref{main theorem}} 
\begin{lemma1}
    \label{mountain pass geometry}
    The functional $I_{\mu}$ has a mountain pass geometry.
\end{lemma1}

\begin{proof}
    From (\ref{1,estimate of the logarithmic nonlinearity}), Lemma \ref{estimate of the exponential nonlinearity} (ii),
    the H$\ddot{\text{o}}$lder inequality, Theorem \ref{fractional Trudinger-Moser inequality}
    and Sobolev embedding inequality, we have that 
    \begin{equation*}
        \begin{split}
            I_{\mu}(u)
            & = \frac{s}{N}\|u\|^{N/s} + \frac{2}{q^2} \int_{\Omega} |u|^q \,dx 
            - \frac{1}{q} \int_{\Omega} |u|^q \ln |u|^2 \,dx - \mu \int_{\Omega} F(u) \,dx\\
            & \geq \frac{s}{N} \|u\|^{N/s} - \frac{\epsilon}{q} \int_{\Omega} |u|^{\theta_1}\,dx
            - \frac{C_{\epsilon}}{q} \int_{\Omega} |u|^{\theta_2} \,dx - \frac{\mu s}{N \sigma} \int_{\Omega} |u|^{N \sigma/s} \,dx\\
            &- \mu \int_{\Omega} |u|^{\theta} e^{|u|^{N/(N-s)}} \,dx\\
            & \geq \frac{s}{N} \|u\|^{N/s} - \frac{\epsilon}{q} \int_{\Omega} |u|^{\theta_1}\,dx
            - \frac{C_{\epsilon}}{q} \int_{\Omega} |u|^{\theta_2} \,dx - \frac{\mu s}{N \sigma} \int_{\Omega} |u|^{N \sigma/s} \,dx \\
            & - (\int_{\Omega} |u|^{t' \theta} \,dx)^{1/t'} (\int_{\Omega} e^{t|u|^{N/(N-s)}} \,dx)^{1/t}\\
            & \geq \frac{s}{N} \|u\|^{N/s} - \frac{\epsilon}{q} C_1 \|u\|^{\theta_1} - \frac{C_{\epsilon}}{q} C_2 \|u\|^{\theta_2}
            - \frac{\mu s}{N \sigma} C_3 \|u\|^{N \sigma/s} - C_4 \|u\|^{\theta},
        \end{split}
    \end{equation*}
    where $C_i, i=1,2,3,4$ are some positive constants and $t,t'>1$ with $\frac{1}{t}+\frac{1}{t'}=1$.
    Since $q,\theta_1, \theta_2, \theta>\frac{N}{s}$, there exist $\rho, \alpha>0$ such that $I_{\mu}(u) \geq \alpha$ for $\|u\|=\rho$.

    Let $u_0 \in W_0^{s,N/s}(\Omega)$ with $\|u_0\| = 1$. Thus it follows from Lemma \ref{estimate of the exponential nonlinearity} (iii)
    and (\ref{2,estimate of the logarithmic nonlinearity}).
    \begin{equation*}
        \begin{split}
            I_{\mu}(tu_0) 
            & = \frac{s}{N}\|t u_0\|^{N/s} + \frac{1}{q^2} \int_{\Omega} [2|tu_0|^q - q|tu_0|^q \ln |tu_0|^2] \,dx
            - \mu \int_{\Omega} F(t u_0) \,dx\\
            & \leq \frac{s}{N}\|t u_0\|^{N/s} + \frac{2 |\Omega|}{q^2}
             - \frac{\mu s}{N \sigma}\int_{\Omega} |t u_0|^{N\sigma/s} \,dx - \frac{\mu}{\theta} \int_{\Omega} |t u_0|^{\theta} \,dx\\
            & \to -\infty,\,\, \text{as}\, t \to \infty.
        \end{split}
    \end{equation*}
    Then we deduce that $I_{\mu}(t_0 u_0)<0$ and $\|t_0 u_0 \| > \rho$ for $t_0>0$ large enough, which completes the proof.
\end{proof}

\begin{proof}[Proof of Theorem \ref{main theorem}]
    Let 
    \begin{equation*}
        \Gamma = \{ \gamma \in C([0,1], W_0^{s,N/s}(\Omega)): \gamma(0)=0 \,\,\text{and}\,\, \gamma(1)=t_0 u_0 \}
    \end{equation*}
    be the class of joining $0$ and $t_0 u_0$.
    Next, we claim that 
    \begin{equation}
        \label{the range of c_mu}
        c_{\mu} \coloneqq \inf_{\gamma \in \Gamma} \max_{t\in [0,1]} I_{\mu}(\gamma(t))
        < (\frac{s}{N} - \frac{1}{q}) \alpha_{s,N}^{(N-s)/s},
    \end{equation}
    for large $\mu$.
    Assuming that (\ref{the range of c_mu}) holds true, we can deduce from Lemma \ref{mountain pass geometry},
    Lemma \ref{PS_c condition} and the Mountain Pass Theorem that there exists at least one nontrivial critical point of $I_{\mu}$.

    Note that $\lim_{t\to \infty}I_{\mu}(tu_0) = -\infty$ for $u_0 \in W_0^{s,N/s}(\Omega)$ with $\|u_0\| = 1$, from which
    there holds $\sup_{t\geq0} I_{\mu}(t u_0) = I_{\mu}(t_{\mu}u_0)$ for some $t_{\mu}>0$. Therefore, 
    by Lemma \ref{estimate of the exponential nonlinearity} (i), (iii),
     $t_{\mu}$ satisfies that 
    \begin{equation}
        \label{t_mu}
        \begin{split}
        t_{\mu}^{N/s} 
        & = \int_{\Omega} |t_{\mu}u_0|^q \ln |t_{\mu}u_0|^2 \,dx + \mu \int_{\Omega} f(t_{\mu}u_0)t_{\mu}u_0 \,dx\\
        & \geq \int_{\Omega} |t_{\mu}u_0|^q \ln |t_{\mu}u_0|^2 \,dx + 
        \frac{sq}{N \sigma} \mu \int_{\Omega} |t_{\mu}u_0|^{N \sigma/s} \,dx +
        \frac{q}{\theta} \mu \int_{\Omega} |t_{\mu}u_0|^{\theta} \,dx.
        \end{split}
    \end{equation}
    The sequence $(t_\mu)$ is bounded since $\theta > N/s$. Now we claim that $t_{\mu} \to 0$ as $\mu \to \infty$.
    Suppose, by contradiction, that there exist $t_{\mu_0}>0$ and a sequence $\mu_n \to \infty$ as $n \to \infty$
    such that $t_{\mu_n} \to t_{\mu_0}$ as $n \to \infty$. By the Lebesgue's convergence theorem and (\ref{1,estimate of the logarithmic nonlinearity}),
    we have that $\int_{\Omega} |t_{\mu_n}u_0|^q \ln |t_{\mu_n}u_0|^2 \,dx \to \int_{\Omega} |t_{\mu_0}u_0|^q \ln |t_{\mu_0}u_0|^2 \,dx$
    as $n \to \infty$. It follows that $t_{\mu_0}^{N/s} \geq \infty$, which is absurd. Therefore, $t_{\mu} \to 0$ as $\mu \to \infty$.

    We deduce from (\ref{t_mu}) that $\lim_{\mu \to \infty} \mu \int_{\Omega}f(t_{\mu} u_0)t_{\mu} u_0 \,dx=0$
    and 
    \begin{equation*}
        \lim_{\mu \to \infty} [\sup_{t \geq 0}I_{\mu}(t u_0)] = \lim_{\mu \to \infty} I_{\mu}(t_{\mu}u_0)=0.
    \end{equation*}
    Then there exists $\mu_*>0$ such that $c_{\mu} \leq \sup_{t \geq 0}I_{\mu}(t u_0) < (\frac{s}{N} - \frac{1}{q}) \alpha_{s,N}^{(N-s)/s}$.
    The proof of Theorem \ref{main theorem} is now complete.
\end{proof} 
\bibliographystyle{abbrv}
\bibliography{template.bib}

\end{spacing}
\end{document}